\documentclass[12pt]{article}
\usepackage{amsmath, amsthm, amssymb}
\numberwithin{equation}{section}
\begin{document}
\author{Ajai Choudhry}
\title{Brahmagupta quadrilaterals with equal perimeters and equal areas}
\date{}
\maketitle

\begin{abstract}
A cyclic quadrilateral is called a Brahmagupta quadrilateral if its four sides, the two diagonals and the area are all given by integers.
In this paper we consider the hitherto unsolved problem of finding two Brahmagupta quadrilaterals with equal perimeters and equal areas. 
We obtain  two parametric solutions of the problem  -- the first solution generates   examples in which each quadrilateral  has two equal sides while  the second solution gives quadrilaterals all of whose sides are unequal. We also show how more parametric solutions of the problem may be obtained.

\end{abstract}



\section{Introduction}

The problem of finding two or more Heron triangles with equal perimeters and equal areas has been considered  by several mathematicians (\cite{Aa}, \cite{Ch} \cite[p. 295]{Gu2}, \cite{Yi}). It would be recalled that a Heron triangle is defined as a triangle having integer sides and integer area. In this paper we consider  an  analogous problem concerning Brahmagupta quadrilaterals. A cyclic quadrilateral is called a Brahmagupta quadrilateral if its four sides, the two diagonals and the area are all given by integers \cite{Sas}. Till now this problem has not been considered in the literature.

If $a_1,\,a_2,\,a_3,\,a_4$ are the consecutive sides of a cyclic quadrilateral, the formulae for the area $K$ and the two diagonals $d_1$ and $d_2$, given by Brahmagupta \cite[p. 187]{Ev} in the seventh century A.D., are as follows:
\begin{equation}
\begin{aligned}
K&=\sqrt{(s-a_1)(s-a_2)(s-a_3)(s-a_4)},\\
d_1&=\sqrt{(a_1a_2+a_3a_4)(a_1a_3+a_2a_4)/(a_1a_4+a_2a_3)},\\
d_2&=\sqrt{(a_1a_3+a_2a_4)(a_1a_4+a_2a_3)/(a_1a_2+a_3a_4)}, 
\end{aligned}
\label{Brahmaformulae}
\end{equation}
where $s$ is the semi-perimeter, that is,
\begin{equation}
s=(a_1+a_2+a_3+a_4)/2. \label{vals}
\end{equation}
Further,  the radius $R$ of the circumcircle of the cyclic quadrilateral is given by the following formula of Paramesvara \cite{Gu1}:
\begin{equation}
R=\frac{\displaystyle 1}{\displaystyle 4}\ \sqrt{\frac{\displaystyle(a_1a_2+a_3a_4)(a_1a_3+a_2a_4)/(a_1a_4+a_2a_3)}{\displaystyle (s-a_1)(s-a_2)(s-a_3)(s-a_4)}} .
\end{equation}

In this paper we find two parametric solutions of the problem of finding two Brahmagupta quadrilaterals with equal perimeters and equal areas. We note that given an arbitrary Brahmagupta quadrilateral, we may change the order of the sides to get two other Brahmagupta quadrilaterals which will have the same perimeter and same area as the given Brahmagupta quadrilateral. The problem is therefore really to find two Brahmagupta quadrilaterals with the stipulated properties and such that the sides of the two quadrilaterals  are given by distinct quadruples of integers.   In Section 2.1 we obtain the first parametric solution that generates   examples in which each quadrilateral  has two equal sides while in Section 2.2, we obtain the second solution which gives quadrilaterals all of whose sides are unequal. In addition to the sides, the diagonals and the areas, the circumradii of the quadrilaterals are also given by integers. 

We note that if there exists a cyclic quadrilateral  whose sides, diagonals and the area are given by rational numbers, by appropriate scaling we can readily obtain  a cyclic quadrilateral whose  sides, diagonals and the area are given by integers. We will therefore refer to such cyclic quadrilaterals also as Brahmagupta quadrilaterals.

\section{Two Brahmagupta quadrilaterals  with \\equal perimeters and equal areas}
 In view of \eqref{Brahmaformulae}, the area $K$ of a cyclic quadrilateral may be written explicitly in terms of the four sides as follows:
\begin{multline}
K= \{(a_1-a_2+a_3+a_4)(a_1+a_2-a_3+a_4)\\
\times (a_1+a_2+a_3-a_4)(-a_1+a_2+a_3+a_4)/16\}^{1/2}. \label{area1}
\end{multline}

Let  $a_1,\,a_2,\,a_3,\,a_4$ and  $b_1,\,b_2,\,b_3,\,b_4$ be the consecutive sides of two Brahmagupta quadrilaterals with  equal perimeters and equal areas. It follows from the formulae \eqref{Brahmaformulae}, \eqref{vals} and \eqref{area1} that $a_i,\;b_i$ must satisfy the following conditions:
\begin{equation}
a_1+a_2+a_3+a_4=b_1+b_2+b_3+b_4, \label{eqside}
\end{equation}
\begin{multline}
(a_1-a_2+a_3+a_4)(a_1+a_2-a_3+a_4)\\
\times(a_1+a_2+a_3-a_4)(-a_1+a_2+a_3+a_4)\\
 =(b_1-b_2+b_3+b_4)(b_1+b_2-b_3+b_4)(b_1+b_2+b_3-b_4)(-b_1+b_2+b_3+b_4) =16K^2, \label{eqarea}
\end{multline}
\begin{equation}
\begin{aligned}
(a_1a_2+a_3a_4)(a_1a_3+a_2a_4)(a_1a_4+a_2a_3)&=u^2,\\
(b_1b_2+b_3b_4)(b_1b_3+b_2b_4)(b_1b_4+b_2b_3)&=v^2,
\end{aligned}
\label{eqdiag}
\end{equation}
where $a_i,\,b_i,\,u,\,v$ and $K$ are all positive integers.

In addition to the above conditions, for the construction of  a cyclic quadrilateral with sides $a_i$,  the values of $a_i,\;i=1,\,\ldots,\,4$, must be such that the sum of any three of the numbers $a_i$ must be greater than the fourth one, and similarly the sum of any three of the four numbers $b_i$ must be greater than the fourth one \cite[Theorem 3.21, p. 57]{CG}. Further, we note that if we obtain a  solution in rational numbers of the diophantine equations \eqref{eqside}, \eqref{eqarea} and \eqref{eqdiag}, by appropriate scaling we can obtain a solution in integers. 

To solve the simultaneous diophantine equations \eqref{eqside}, \eqref{eqarea} and \eqref{eqdiag}, we first make the invertible linear transformation defined by the following relations:
\begin{equation}
\begin{aligned}
a_1&= (-x_1+x_2+x_3+x_4)/2, &a_2& = (x_1-x_2+x_3+x_4)/2,\\
 a_3 &= (x_1+x_2-x_3+x_4)/2, &a_4& = (x_1+x_2+x_3-x_4)/2,\\
b_1&= (-y_1+y_2+y_3+y_4)/2, &b_2& = (y_1-y_2+y_3+y_4)/2,\\
 b_3 &= (y_1+y_2-y_3+y_4)/2, &b_4& = (y_1+y_2+y_3-y_4)/2.
\end{aligned}
\label{ltt1}
\end{equation}

The equations \eqref{eqside}, \eqref{eqarea} and \eqref{eqdiag} may now be written as follows:
\begin{align}
x_1+x_2+x_3+x_4&=y_1+y_2+y_3+y_4, \label{eqsidex}\\
x_1x_2x_3x_4 &= y_1y_2y_3y_4=K^2,\label{eqareax}
\end{align}
\begin{align}
(x_1x_2+x_3x_4)(x_1x_3+x_2x_4)(x_1x_4+x_2x_3)&=u^2, \label{eqdiag1x}\\
(y_1y_2+y_3y_4)(y_1y_3+y_2y_4)(y_1y_4+y_2y_3)&=v^2. \label{eqdiag2x}
\end{align}

\subsection{}
We will now obtain pairs of Brahmagupta quadrilaterals with equal perimeters and equal areas such that each  quadrilateral has two sides equal. It follows from \eqref{ltt1} that two sides of each quadrilateral wll be equal  if and only if two of the numbers $\{x_i\}$ as well as two of the numbers $\{y_i\}$ are equal. There is no loss of generality in taking $x_4=x_3$ and $y_4=y_3$. 

On writing $x_4=x_3$ and $y_4=y_3$, Eq.~\eqref{eqareax} reduces to $x_1x_2x_3^2=y_1y_2y_3^2=K^2$. It follows that both $x_1x_2$ and $y_1y_2$ must be perfect squares, and so we write
\begin{equation}
x_1=p_1q_1^2,\quad x_2=p_1q_2^2,\quad y_1=p_2r_1^2,\quad  y_2=p_2r_2^2, \label{tryx}
\end{equation}
where $p_i,\,q_i,\,r_i,\;i=1,\,2$ are arbitrary parameters. It now follows from Eq.~\eqref{eqareax} that $(p_1q_1q_2x_3)^2 =( p_2r_1r_2y_3)^2$ and we accordingly take
\begin{equation}
x_3=p_2r_1r_2,\quad y_3=p_1q_1q_2. \label{tryx2}
\end{equation}

Further, with the values of $x_i,\,y_i$ given by \eqref{tryx} and \eqref{tryx2}, Eq.~\eqref{eqsidex} reduces to
\begin{equation}
(q_1-q_2)^2p_1-(r_1-r_2)^2p_2=0,
\end{equation}
and accordingly, we take
\begin{equation}
p_1= (r_1-r_2)^2,\quad \quad p_2=(q_1-q_2)^2. \label{valp12}
\end{equation}

With the values of $x_i,\,y_i,\,p_i$ given by \eqref{tryx}, \eqref{tryx2} and \eqref{valp12}, the conditions \eqref{eqsidex} and \eqref{eqareax} are satisfied, while  Eqs.~\eqref{eqdiag1x} and \eqref{eqdiag2x} reduce to the following two conditions respectively:
\begin{align}
(q_1-q_2)^4(q_1^2+q_2^2)^2r_1^2r_2^2(r_1-r_2)^4\phi(q_1,\,q_2,\,r_1,\,r_2) &= u^2, \label{eqdiag1xnew}\\
q_1^2q_2^2(q_1-q_2)^4(r_1-r_2)^4(r_1^2+r_2^2)^2 \phi(q_1,\,q_2,\,r_1,\,r_2) &= v^2, \label{eqdiag2xnew}
\end{align}
where
\begin{multline}
\phi(q_1,\,q_2,\,r_1,\,r_2)=r_1^2r_2^2q_1^4-4r_1^2r_2^2q_1^3q_2+(r_1^4-4r_1^3r_2+12r_1^2r_2^2\\
-4r_1r_2^3+r_2^4)q_1^2q_2^2-4r_1^2r_2^2q_1q_2^3+r_1^2r_2^2q_2^4.\quad \quad \quad 
\end{multline}

It  follows that if we choose $q_1$ and $q_2$ such that $\phi(q_1,\,q_2,\,r_1,\,r_2)$ becomes a perfect square, then  both conditions the \eqref{eqdiag1xnew} and \eqref{eqdiag2xnew} will be satisfied. Now  $\phi(q_1,\,q_2,\,r_1,\,r_2)$ is a quartic function of $q_1$ and $q_2$ and we can readily find suitable  values of $q_1,\,q_2$ by following a method described by Fermat (as quoted by Dickson \cite[p. 639]{Di}).  We thus obtain    the following values of $q_1,\,q_2$:
\begin{equation}
q_1=8r_1^2r_2^2,\quad q_2=r_1^4-4r_1^3r_2+10r_1^2r_2^2-4r_1r_2^3+r_2^4. \label{valq}
\end{equation}

With the above   values of  $q_1,\,q_2$, the relations \eqref{valp12} yield the values of $p_1,\,p_2$ and now using  the relations \eqref{tryx} and \eqref{tryx2}, we obtain the following values of $x_i,\,y_i$. Finally, using the relations \eqref{ltt1}, we get the sides $a_i,\,b_i,\;i=1,\,\ldots,\,4$, of two Brahmagupta quadrilaterals with equal perimeters and equal areas. These values of $a_i,\,b_i,\;i=1,\,\ldots,\,4$, are given below in terms of two arbitrary parameters $r_1$ and $r_2$:
\begin{equation}
\begin{aligned}
a_1 &= (r_1^2+r_2^2)(r_1^2-4r_1r_2+r_2^2)(r_1^6-4r_1^5r_2+\\
& \quad \quad 19r_1^4r_2^2-40r_1^3r_2^3+19r_1^2r_2^4-4r_1r_2^5+r_2^6), \\
a_2 &= -(r_1^2+r_2^2)(r_1^2-4r_1r_2+r_2^2)(r_1^6-8r_1^5r_2\\
& \quad \quad+35r_1^4r_2^2-48r_1^3r_2^3+35r_1^2r_2^4-8r_1r_2^5+r_2^6), \\
a_3 = a_4&=(r_1-r_2)^2(r_1^8-8r_1^7r_2+36r_1^6r_2^2-88r_1^5r_2^3\\
& \quad \quad +198r_1^4r_2^4-88r_1^3r_2^5+36r_1^2r_2^6 -8r_1r_2^7+r_2^8), \\
b_1& = -(r_1-r_2)(r_1^4-8r_1^3r_2+10r_1^2r_2^2+r_2^4)\\
& \quad \quad \times(r_1^5+r_1^4r_2-6r_1^3r_2^2+18r_1^2r_2^3-7r_1r_2^4+r_2^5),\\
 b_2 &= (r_1-r_2)(r_1^4+10r_1^2r_2^2-8r_1r_2^3+r_2^4)\\
& \quad \quad \times(r_1^5-7r_1^4r_2+18r_1^3r_2^2-6r_1^2r_2^3+r_1r_2^4+r_2^5),\\
 b_3 = b_4&=(r_1^2-4r_1r_2+r_2^2)^2(r_1^2+r_2^2)^3.
\end{aligned}
\label{valab}
\end{equation}

The common perimeter of the two Bramhagupta quadrilaterals obtained above is given by 
\begin{multline}
2r_1^{10}-16r_1^9r_2+74r_1^8r_2^2-256r_1^7r_2^3+724r_1^6r_2^4-992r_1^5r_2^5+724r_1^4r_2^6\\
-256r_1^3r_2^7+74r_1^2r_2^8-16r_1r_2^9+2r_2^{10},\quad \quad \quad \quad 
\end{multline}
while the common area is given by 
\begin{multline}
\quad \quad \quad 32r_1^3r_2^3(r_1-r_2)^2(r_1^2+r_2^2)^2(r_1^2-4r_1r_2+r_2^2)^2\\
\times (r_1^4-4r_1^3r_2+10r_1^2r_2^2-4r_1r_2^3+r_2^4). \quad \quad 
\end{multline}

The two    diagonals of  the first  quadrilateral are of lengths
\begin{multline}
2r_1r_2(r_1^8-8r_1^7r_2+52r_1^6r_2^2-152r_1^5r_2^3+230r_1^4r_2^4-152r_1^3r_2^5+52r_1^2r_2^6-8r_1r_2^7+r_2^8). \label{diag1}
\end{multline}
and 
\begin{multline}
2(r_1^8-8r_1^7r_2+36r_1^6r_2^2-88r_1^5r_2^3+198r_1^4r_2^4-88r_1^3r_2^5+36r_1^2r_2^6-8r_1r_2^7+r_2^8)\\
\quad \times (r_1-r_2)^2(r_1^2+r_2^2)^2(r_1^2-4r_1r_2+r_2^2)^2/(r_1^8-8r_1^7r_2+52r_1^6r_2^2\\
\quad -152r_1^5r_2^3+230r_1^4r_2^4-152r_1^3r_2^5+52r_1^2r_2^6-8r_1r_2^7+r_2^8),
\end{multline}
while the two    diagonals of  the second  quadrilateral have  lengths given by \eqref{diag1} and
\begin{multline}
16r_1r_2(r_1^4-4r_1^3r_2+10r_1^2r_2^2-4r_1r_2^3+r_2^4)(r_1-r_2)^2(r_1^2-4r_1r_2+r_2^2)^2(r_1^2+r_2^2)^3\\
\quad \times (r_1^8-8r_1^7r_2+52r_1^6r_2^2-152r_1^5r_2^3+230r_1^4r_2^4-152r_1^3r_2^5+52r_1^2r_2^6-8r_1r_2^7+r_2^8)^{-1}.
\end{multline}

The circumradii of the two quadilaterals are given by 
\begin{multline}
(r_1^8-8r_1^7r_2+36r_1^6r_2^2-88r_1^5r_2^3+198r_1^4r_2^4-88r_1^3r_2^5+36r_1^2r_2^6-8r_1r_2^7+r_2^8)\\
\times (r_1^8-8r_1^7r_2+52r_1^6r_2^2-152r_1^5r_2^3+230r_1^4r_2^4-152r_1^3r_2^5+52r_1^2r_2^6-8r_1r_2^7+r_2^8)\\
\times \{16r_1r_2(r_1^4-4r_1^3r_2+10r_1^2r_2^2-4r_1r_2^3+r_2^4)\}^{-1}
\end{multline}
and
\begin{multline}
(r_1^2+r_2^2)(r_1^8-8r_1^7r_2+52r_1^6r_2^2-152r_1^5r_2^3+\\
230r_1^4r_2^4-152r_1^3r_2^5+52r_1^2r_2^6-8r_1r_2^7+r_2^8)/2,
\end{multline}
respectively.

It is readily seen that when we choose $r_1,\,r_2$ such that $1.63 < r_1/r_2 < 2.11$, the values of $a_i,\,b_i$ given by \eqref{valab} are all positive, the sum any three of the numbers $\{a_i,\;i=1,\ldots,\,4\}$ is greater than the fourth one and the same is true of the quadruple $\{b_i,\;i=1,\ldots,\,4\}$, and thus there actually exist quadrilaterals whose sides are given by the quadruples $\{a_i,\;i=1,\ldots,\,4\}$ and $\{b_i,\;i=1,\ldots,\,4\}$.

As a numerical example, when $r_1=2,\,r_2=1$, we get two Brahmagupta quadrilaterals with the following sides:
\begin{equation}
\begin{aligned}
a_1 &= 165, \quad a_2 &= 1635, \quad a_3 &=a_4&= 1313,\\
 b_1 &= 413, \quad b_2 &= 1763, \quad b_3 &= b_4&=1125.
\end{aligned}
\end{equation}
 The perimeters and areas of both of these quadrilaterals are given by $4426$ and  $979200$ respectively.

The diagonals of the first  quadrilateral have lengths 1412 and 590850/353 while the lengths of the diagonals of the second quadrilateral are 1412 and 612000/353. The circumradii of the two quadrilaterals are 
463489/544 and 1765/2 respectively. By appropriate scaling, we can readily get two Brahmagupta quadrilaterals  whose diagonals and circumradii are also integers.

The lengths of the  diagonals of the two quadrilaterals given above have been computed on the basis that the consecutive sides of the two quadrilaterals are $a_1,\,a_2,\,a_3,\,a_4$ and  $b_1,\,b_2,\,b_3,\,b_4$ as was stated initially. If we take the consecutive sides of the two quadrilaterals as $a_1,\,a_3,\,a_2,\,a_4$ and  $b_1,\,b_3,\,b_2,\,b_4$ respectively,  we get two trapeziums with the same perimeter and same  area as the first pair of quadrilaterals while the  four diagonals of the two trapeziums are all of the same length which is given by \eqref{diag1}.

We note that while we have obtained one parametric solution using the values of $q_1,\,q_2$ given by \eqref{valq}, by repeatedly applying the aforementioned method of Fermat, we can get infinitely many pairs of values of $q_1,\,q_2$ such that $\phi(q_1,\,q_2,\,r_1,\,r_2)$ becomes a perfect square, and thus we can obtain infinitely many parametric solutions giving pairs of Brahmagupta quadrilaterals with equal perimeters and equal areas and such that each quadrilateral has two sides equal.

\subsection{} We will now find a solution of the simultaneous diophantine equations  \eqref{eqsidex}, \eqref{eqareax}, \eqref{eqdiag1x} and \eqref{eqdiag2x} in which the quadruples $\{x_i,\;i=1,\ldots,4\}$ and $\{y_i,\;i=1,\ldots,4\}$ consist of distinct integers.

To solve Eqs.~\eqref{eqsidex} and \eqref{eqareax}, we first obtain two quadruples of integers with equal sums and equal products by applying a method described by Choudhry \cite[Lemma 5, p. 769]{Ch3}. Accordingly, we solve the equation,
\begin{equation}
(p_1-p_2)/(q_1-q_2) = (r_1-r_2)/(s_1-s_2), \label{aux1}
\end{equation}
and take the integers $x_i,\,y_i,\;i=1,\ldots,4$, as follows:
\begin{equation}
\begin{aligned} 
x_1&=p_1s_1, &x_2&=p_2s_2, &x_3&=q_1r_2, &x_4&=q_2r_1,\\
y_1&=p_1s_2, &y_2&=p_2s_1, &y_3&=q_1r_1, &y_4&=q_2r_2. 
\end{aligned}
\label{valxy}
\end{equation}
A solution of \eqref{aux1} is given by
\begin{equation}
q_2 = -mp_1+mp_2+q_1, s_2 = -mr_1+mr_2+s_1,
\end{equation}
where $m$ is an arbitrary parameter, and we thus get the two quadruples $\{x_i\}$ and $\{y_i\}$ given by
\begin{equation}
\begin{aligned}
x_1 &= p_1s_1, &x_2 &= -p_2(mr_1-mr_2-s_1), \\
x_3 &= r_2q_1, &x_4 &= -r_1(mp_1-mp_2-q_1), \\
y_1 &= -p_1(mr_1-mr_2-s_1), &y_2 &= p_2s_1,\\
 y_3 &= r_1q_1, &y_4 &= -r_2(mp_1-mp_2-q_1)
\end{aligned}
\label{val2xy}
\end{equation}
where $m,\,p_1,\,p_2,\,q_1,\,r_1,\,r_2$ and $s_1$ are arbitrary parameters. 

With the values of $x_i,\,y_i$, given by \eqref{val2xy}, Eq.~\eqref{eqsidex} is satisfied while Eq.~\eqref{eqareax} reduces to the equation,
\begin{equation}
p_1p_2q_1r_1r_2s_1(mp_1-mp_2-q_1)(mr_1-mr_2-s_1)=K^2. \label{eqareax2}
\end{equation}
We now choose
\begin{equation}
q_1 = n(p_1-p_2),\quad  s_1 = n(r_1-r_2), \label{valq1s1}
\end{equation}
when Eq.~\eqref{eqareax2} reduces to the condition,
\begin{equation}
p_1p_2r_1r_2(m-n)^2n^2(p_1-p_2)^2(r_1-r_2)^2=K^2. \label{eqareax3}
\end{equation}
A solution of this equation is readily obtained and may be written as,
\begin{equation}
p_1=g^2r_1,\quad p_2=h^2r_2,\quad K=ghr_1r_2(m-n)n(r_1-r_2)(g^2r_1-h^2r_2), \label{valr1r2t}
\end{equation}
where $g$ and $h$ are arbitrary nonzero parameters. 

We have now obtained a solution of   Eqs.~\eqref{eqsidex} and \eqref{eqareax}, and using the relations \eqref{val2xy}, \eqref{valq1s1} and \eqref{valr1r2t}, the values of $x_i,\,y_i$  may be written explicitly   as follows:
\begin{equation}
\begin{aligned}
x_1 &= g^2r_1n(r_1-r_2), &x_2 &= -h^2r_2(r_1-r_2)(m-n),\\
 x_3 &= r_2n(g^2r_1-h^2r_2), &x_4 &= -r_1(g^2r_1-h^2r_2)(m-n),\\
 y_1 &= -g^2r_1(r_1-r_2)(m-n), &y_2 &= h^2r_2n(r_1-r_2), \\
y_3 &= r_1n(g^2r_1-h^2r_2), &y_4 &= -r_2(g^2r_1-h^2r_2)(m-n),
\end{aligned}
\label{solxy}
\end{equation}
where $g,\,h,\,m,\,n,\,r_1$ and $r_2$ are arbitrary parameters.

With the values of $x_i,\,y_i$, defined by \eqref{solxy},  Eqs.~\eqref{eqdiag1x} and \eqref{eqdiag2x} reduce to the following two conditions respectively:
\begin{multline}
(m-n)^2n^2r_1^2r_2^2(r_1-r_2)^2(g^2r_1-h^2r_2)^2(g^2r_1^2+h^2r_2^2)
(g^4r_1^2\\+g^2h^2r_1^2-4g^2h^2r_1r_2+g^2h^2r_2^2+h^4r_2^2)\{(g^2n^2+h^2(m-n)^2\}=u^2, \label{eqdiag1x2}
\end{multline}
and
\begin{multline}
(m-n)^2n^2r_1^2r_2^2(r_1-r_2)^2(g^2r_1-h^2r_2)^2(g^2r_1^2+h^2r_2^2)
(g^4r_1^2\\ \quad +g^2h^2r_1^2-4g^2h^2r_1r_2+g^2h^2r_2^2+h^4r_2^2)\{g^2(m-n)^2+h^2n^2\}=v^2. \label{eqdiag2x2}
\end{multline}
On multiplying Eqs.~\eqref{eqdiag1x2} and \eqref{eqdiag2x2} by $v^2$ and $u^2$ respectively, and taking the difference, we get, after removing certain factors, the condition,
\begin{equation}
(mu-nu-nv)(mu-nu+nv)g^2-(mv+nu-nv)(mv-nu-nv)h^2 = 0, \label{cond1}
\end{equation}
which further reduces, on writing
\begin{equation}
m=n(1-X),\quad g=hY/((Xu-v)(Xu+v),
\end{equation}
to the condition,
\begin{equation}
Y^2=(Xu-v)(Xu+v)(Xv-u)(Xv+u). \label{ecXY}
\end{equation}

For rational values of $u$ and $v$, Eq.~\eqref{ecXY} represents a quartic model of an elliptic curve which, in general, is of rank 0. With some experimentation, we found that when we take 
\begin{equation}
u = t(t^4-2t^2+5), \quad v=5t^4-2t^2+1,
\end{equation}
where $t$ is an arbitrary rational parameter, the elliptic curve \eqref{ecXY} is, in general, of positive rank, and a point of infinite order on the curve \eqref{ecXY} is given by
\begin{equation}
\begin{aligned}
X&=(3t^4-6t^2-1)/\{t(t^4+6t^2-3)\},\\
Y&=4(t-1)(t+1)(t^2+1)^2(t^2+2t-1)(t^2-2t-1)\\
& \quad \quad \times (t^8+20t^6-26t^4+20t^2+1)/\{t(t^4+6t^2-3)^2\}.
\end{aligned}
\end{equation}
We thus obtain the following solution of Eq.~\eqref{cond1}:
\begin{equation}
\begin{aligned}
g &= (t^2+2t-1)(t^2-2t-1), &h &= -4t(t^2+1), \\
m &= t^5-3t^4+6t^3+6t^2-3t+1, &n &= (t^4+6t^2-3)t, \\
u &= t(t^4-2t^2+5)w, &v &= (5t^4-2t^2+1)w,
\end{aligned}
\label{valgh}
\end{equation}
where $t$ and $w$ are arbitrary parameters.

Substituting the values of $g,\,h,\,m,\,n$ and $u$ in Eq.~\eqref{eqdiag1x2}, we get the condition,
\begin{equation}
w^2=\phi_1^2(r_1,\,r_2,\,t)\phi_2(r_1,\,r_2,\,t), \label{eqdiag1x3}
\end{equation}
where
\begin{equation}
\begin{aligned}
\phi_1(r_1,\,r_2,\,t)&=r_1r_2(r_1-r_2)t(3t^4-6t^2-1)(t^4+6t^2-3) \\
& \quad \quad \times \{r_1(t^2+2t-1)^2(t^2-2t-1)^2-16r_2t^2(t^2+1)^2\},\\
\phi_2(r_1,\,r_2,\,t)&=(t^4-4t^3+10t^2-4t+1)^2(t^4+4t^3+10t^2+4t+1)^2\\
& \quad \;\; \times (t^2+2t-1)^4(t^2-2t-1)^4r_1^4-64t^2(t^2+1)^2(t^2+2t-1)^4\\
& \quad \;\; \times(t^4-4t^3+10t^2-4t+1)(t^4+4t^3+10t^2+4t+1)\\
& \quad \;\; \times  (t^2-2t-1)^4r_1^3r_2+32t^2(t^2+1)^2(t^2+2t-1)^2\\
&\quad \;\; \times (t^4-4t^3+10t^2-4t+1)^2(t^4+4t^3+10t^2+4t+1)^2\\
& \quad \;\; \times (t^2-2t-1)^2r_1^2r_2^2-1024t^4(t^2+2t-1)^2(t^2-2t-1)^2\\
& \quad \;\; \times (t^4-4t^3+10t^2-4t+1)(t^4+4t^3+10t^2+4t+1)\\
& \quad \;\; \times (t^2+1)^4r_1r_2^3+256t^4(t^4-4t^3+10t^2-4t+1)^2\\
& \quad \;\; \times (t^4+4t^3+10t^2+4t+1)^2(t^2+1)^4r_2^4.
\end{aligned}
\end{equation}
Since $\phi_2(r_1,\,r_2,\,t)$ is a quartic function of $r_1$ and $r_2$, we can readily find   values of $r_1,\,r_2$ that make $\phi_2(r_1,\,r_2,\,t)$ a perfect square by following the  method of Fermat already mentioned in Section 2.1 .  We thus obtain    the following values of $r_1,\,r_2$  that satisfy the condition \eqref{eqdiag1x3}:
\begin{equation}
\begin{aligned}\
r_1&=32t^2(t^4-4t^3+10t^2-4t+1)(t^4+4t^3+10t^2+4t+1)(t^2+1)^2,\\
r_2&=(t^8+4t^7+4t^6-20t^5+70t^4-20t^3+4t^2+4t+1)\\
& \quad \;\; \times(t^8-4t^7+4t^6+20t^5+70t^4+20t^3+4t^2-4t+1).
\end{aligned}
\label{valr}
\end{equation}

With the values of $g,\,h,\,m,\,n$ and $r_1,\,r_2$ given by \eqref{valgh} and \eqref{valr} respectively, we get the values of $x_i,\,y_i,\;i=1,\,\ldots,\,4$, from the relations \eqref{solxy} and we thus obtain a solution of the diophantine system consisting of the equations 
\eqref{eqsidex}, \eqref{eqareax}, \eqref{eqdiag1x} and \eqref{eqdiag2x}. Finally, using the relations \eqref{ltt1}, we obtain, in terms of an arbitrary parameter $t$,  the following values of $a_i,\,b_i,\;i=1,\,\ldots,\,4$, satisfying the conditions \eqref{eqside}, \eqref{eqarea} and \eqref{eqdiag}:
\begin{equation}
\begin{aligned}
a_1 &= (t^4+4t^3+10t^2+4t+1)(t^{16}-24t^{14}+48t^{13}-100t^{12}-672t^{11}\\
&\quad \;\;-1128t^{10}+2960t^9-3002t^8-4032t^7+1240t^6+592t^5-868t^4\\
&\quad \;\;+96t^3+40t^2-16t+1)(3t^{17}-15t^{16}+8t^{15}+136t^{14}-140t^{13}\\
&\quad \;\;-276t^{12}+2488t^{11}-2792t^{10}+1170t^9-42t^8+2488t^7-1800t^6\\
&\quad \;\;-140t^5-52t^4+8t^3-24t^2+3t+1), \\
a_2 &= -(t^4-4t^3+10t^2-4t+1)(t^{17}+t^{16}-72t^{15}-264t^{14}-580t^{13}\\
&\quad \;\;-964t^{12}-2680t^{11}-2232t^{10}-250t^9+1542t^8-2680t^7-952t^6\\
&\quad \;\;-580t^5+60t^4-72t^3-8t^2+t+1)(t^{16}-8t^{14}+60t^{12}-768t^{11}\\
&\quad \;\;-1720t^{10}+1542t^8+2560t^7+328t^6+2048t^5+1084t^4+256t^3\\
&\quad \;\;-8t^2+1),\\
 a_3 &= -(t^4-4t^3+10t^2-4t+1)(t^{16}-24t^{14}-48t^{13}-100t^{12}\\
&\quad \;\;+672t^{11}-1128t^{10}-2960t^9-3002t^8+4032t^7+1240t^6-592t^5\\
&\quad \;\;-868t^4-96t^3+40t^2+16t+1)(3t^{17}+15t^{16}+8t^{15}-136t^{14}\\
&\quad \;\;-140t^{13}+276t^{12}+2488t^{11}+2792t^{10}+1170t^9+42t^8\\
&\quad \;\;+2488t^7+1800t^6-140t^5+52t^4+8t^3+24t^2+3t-1), \\
a_4 &= -(t^4+4t^3+10t^2+4t+1)(t^{16}-8t^{14}+60t^{12}+768t^{11}-1720t^{10}\\
&\quad \;\;+1542t^8-2560t^7+328t^6-2048t^5+1084t^4-256t^3-8t^2+1)\\
&\quad \;\; \times (t^{17}-t^{16}-72t^{15}+264t^{14}-580t^{13}+964t^{12}\\
&\quad \;\;-2680t^{11}+2232t^{10}-250t^9-1542t^8-2680t^7+952t^6\\
&\quad \;\;-580t^5-60t^4-72t^3+8t^2+t-1),
\end{aligned} \label{valaf}
\end{equation}
\begin{equation}
\begin{aligned}
b_1 &= -(t^4+4t^3+10t^2+4t+1)(t^{16}-16t^{15}+40t^{14}+96t^{13}-868t^{12}\\
&\quad \;\;+592t^{11}+1240t^{10}-4032t^9-3002t^8+2960t^7-1128t^6-672t^5\\
&\quad \;\;-100t^4+48t^3-24t^2+1)(t^{17}+3t^{16}-24t^{15}+8t^{14}-52t^{13}\\
&\quad \;\;-140t^{12}-1800t^{11}+2488t^{10}-42t^9+1170t^8-2792t^7+2488t^6\\
&\quad \;\;-276t^5-140t^4+136t^3+8t^2-15t+3), \\
b_2 &= (t^4-4t^3+10t^2-4t+1)(t^{17}+t^{16}-8t^{15}-72t^{14}+60t^{13}-580t^{12}\\
&\quad \;\;-952t^{11}-2680t^{10}+1542t^9-250t^8-2232t^7-2680t^6-964t^5\\
&\quad \;\;-580t^4-264t^3-72t^2+t+1)(t^{16}-8t^{14}+256t^{13}+1084t^{12}+2048t^{11}\\
&\quad \;\;+328t^{10}+2560t^9+1542t^8-1720t^6-768t^5+60t^4-8t^2+1),\\
 b_3 &= -(t^4+4t^3+10t^2+4t+1)(t^{16}-8t^{14}-256t^{13}+1084t^{12}-2048t^{11}\\
&\quad \;\;+328t^{10}-2560t^9+1542t^8-1720t^6+768t^5+60t^4-8t^2+1)\\
&\quad \;\; \times (t^{17}-t^{16}-8t^{15}+72t^{14}+60t^{13}+580t^{12}-952t^{11}+2680t^{10}+1542t^9\\
&\quad \;\;+250t^8-2232t^7+2680t^6-964t^5+580t^4-264t^3+72t^2+t-1),\\
 b_4 &= -(t^4-4t^3+10t^2-4t+1)(t^{16}+16t^{15}+40t^{14}-96t^{13}-868t^{12}\\
&\quad \;\;-592t^{11}+1240t^{10}+4032t^9-3002t^8-2960t^7-1128t^6+672t^5\\
&\quad \;\;-100t^4-48t^3-24t^2+1)(t^{17}-3t^{16}-24t^{15}-8t^{14}-52t^{13}\\
&\quad \;\;+140t^{12}-1800t^{11}-2488t^{10}-42t^9-1170t^8-2792t^7\\
&\quad \;\;-2488t^6-276t^5+140t^4+136t^3-8t^2-15t-3).
\end{aligned} \label{valbf}
\end{equation}

It is readily verified  that when  $ 4.991 < t < 5.565 $, the values of $a_i,\,b_i$ are all positive, the sum any three of the numbers $\{a_i,\;i=1,\ldots,\,4\}$, is greater than the fourth one and the same is true of the quadruple $\{b_i,\;i=1,\ldots,\,4\}$. Thus, we can construct three Brahmagupta quadrilaterals with sides $\{a_i,\;i=1,\ldots,\,4\}$, and three Brahmagupta quadrilaterals with sides $\{b_i,\;i=1,\ldots,\,4\}$, and all of these quadrilaterals will have the same perimeter and the same area. Further, the four sides of these quadrilaterals are unequal, and their circumradii are  rational. 

We do not explicitly give the common perimeter  and the common area of these Brahmagupta quadrilaterals as their values are cumbersome to write. Similarly we omit the lengths of the diagonals and the circumradii of these quadrilaterals. 

As a numerical example, when $t=5$, we get, after removing a common factor, the following values of $a_i,\,b_i:$
\begin{equation}
\begin{aligned}
a_1 &= 1910470516999149312, &a_2 &= 175866555513132912053,\\
 a_3 &= 169314770763852594617, &a_4 &= 207503184245618672382, \\
b_1 &= 154300800756891939924, &b_2 &= 50195745087237056747, \\
b_3 &= 121029496193614687182, &b_4 &= 229068939001859644511.
\end{aligned} \label{valabex2}
\end{equation}
The common perimeter and the common area of both of these quadrilaterals are given by $554594981039603328364$ and  \[14509220341219325824870053111347523537900\] respectively.

Further, the two diagonals and the circumradius of the first quadrilateral are given by 
\[
\begin{aligned}
 &250496054986226007288150003450/1204106621,\\
 &43680775787512057583999745775/246823021,\\
 {\rm and}\quad  &1338548290849915267747645/12376,
\end{aligned}
\]
respectively while the diagonals and the circumradius of the second  quadrilateral are given by
\[
\begin{aligned}
 &123610451156476856682515/769,\\
 &46331747007719685906339040691/246823021,\\
 {\rm and} \quad &7098921625266102351020269/61880,
\end{aligned}\]
respectively.

While the above parametric solution generates infinitely many examples of two Brahmagupta quadrilaterals with equal perimeters and equal areas, we note that by repeated application of the afoementioned method of Fermat, we can obtain infinitely many solutions of Eq.~\eqref{eqdiag1x3} and thus obtain more parametric solutions of the problem. 

\section{An open problem}
It is easy to show that given any positive integer $n$ howsoever large, there exist  $n$ cyclic quadrilaterals with integer sides and having the same perimeter and the same integer area. In view of the computations in Section 2, it suffices to obtain infinitely many solutions of the simultaneous diophantine equations,
\begin{equation}
x_1+x_2+x_3+x_4=k_1,\,\quad \quad x_1x_2x_3x_4=k_2^2, \label{eqk}
\end{equation}
where $k_1$ and $k_2$ are certain constants. 

It has been shown  by Schinzel \cite{Sch} that there exist infinitely many solutions in positive rational numbers of the simultaneous diophantine equations 
\begin{equation}
x_1+x_2+x_3=6,\,\quad \quad x_1x_2x_3=6,
\end{equation}
and taking $x_4=6k^2$ where $k$ is a suitably chosen constant, we immediately get infinitely many solutions in positive rational numbers of the simultaneous diophantine equations \eqref{eqk}, and we thereby get infinitely many cyclic quadrilaterals with rational sides, and having the same perimeter and the same rational area. Appropriate scaling yields  $n$ cyclic quadrilaterals with integer sides and having the required properties. 

It remains an open question whether there for any arbitrary integer $n \geq 3$, there exist $n$ Brahmagupta quadrilaterals  with the same perimeter and the same area and such that the sides of the quadrilaterals are given by distinct quadruples of integers.

\begin{center}
\Large
Acknowledgments
\end{center}
 
I wish to  thank the Harish-Chandra Research Institute, Allahabad for providing me with all necessary facilities that have helped me to pursue my research work in mathematics.

\medskip

\noindent Postal Address: Ajai Choudhry, 
\newline \hspace{1.05 in}
13/4 A Clay Square,
\newline \hspace{1.05 in} Lucknow - 226001, INDIA.
\newline \noindent  E-mail: ajaic203@yahoo.com

\end{document}